\documentclass[11pt]{article}

\usepackage{amsfonts,amsmath,graphicx}
 \graphicspath{{.}}
   \makeatletter
   \let\accent@spacefactor\relax
   \makeatother
    \usepackage[french]{babel}
 \usepackage{amsfonts}
 \usepackage{amsmath}   
 \usepackage{amsthm}    
 \usepackage{amscd}     
 \usepackage{euscript}

 \def\C{{\Bbb C}}
 \def\P{{\Bbb P}}

 \def\Z{{\Bbb Z}}

 \newtheorem{defi}{D\'{e}finition}[section]
 \newtheorem{pro}[defi]{Proposition}
 
 \newtheorem{lem}[defi]{Lemma}
 \newtheorem{theo}[defi]{Th\'eor\`eme}
 
 \newtheorem{coro}[defi]{Corollaire}

\title{Vari\'et\'es de type Togliatti}
\author{Jean Vall\`es}
\date{}

\begin{document}
\maketitle

\section{Introduction } 
La surface de Del Pezzo $S_{6}\subset \P^{6}$ de degr\'e $6$ obtenue par l'\'eclatement de trois points non align\'es du plan projectif complexe $\P^{2}$ (par le syst\`eme lin\'eaire des cubiques) poss\`ede une propri\'et\'e spectaculaire : ses hyperplans osculateurs (espace engendr\'e par les d\'eriv\'ees partielles secondes d'une repr\'esentation param\'etrique)  ont un point commun dans l'espace ambiant. Autrement dit  la surface ``osculatrice'' dans l'espace dual $\P^{\vee 6}$ est d\'eg\'en\'er\'ee. On dit alors qu'elle v\'erifie une \'equation de Laplace. Comme cette propri\'et\'e a \'et\'e d\'emontr\'ee par Togliatti, la surface $S_{6}$ h\'erite du nom de  surface de Togliatti (voir \cite{FI}, \cite{LM} et \cite{T}). 

\smallskip

La surface de Togliatti est une projection de la surface de Veronese $v_{3}(\P^{2})\subset 
\P(H^0(O_{\P^{2}}(3))) $ mais elle peut aussi \^etre d\'efinie comme section hyperplane g\'en\'erale de la vari\'et\'e  de Segre 
$Seg(1,3)\subset \P(H^0(O_{\P^{1}\times \P^{1}\times \P^{1}}(1,1,1)))$.

\smallskip

Dans une premi\`ere partie je reviens rapidement sur ces deux descriptions et je donne une preuve d'une version simplifi\'ee du th\'eor\`eme  de Togliatti (voir th\'eor\`eme \ref{toto}).

\smallskip

Dans les deuxi\`eme et troisi\`eme parties on montre que les espaces $2n$-tangents (voir ci-dessous)  des projections bien choisies des Veronese $v_{2n+1}(\P^{2})\subset 
\P(H^0(O_{\P^{2}}(2n+1))) $ et  des sections hyperplanes g\'en\'erales des Segre 
$Seg(1,2n+1)\subset \P(H^0(O_{\P^{1}\times \cdots \times \P^{1}}(1,\cdots,1)))$
ont un point commun (voir les th\'eor\`emes \ref{vero} et \ref{pun}). Le point cl\'e est
la polarit\'e par rapport aux courbes rationnelles normales et par rapport aux produit de 
$\P^{1}$.
 
\medskip

\textbf{Notations.} On note $T_{x}X\subset \P^N$ l'espace tangent d'une  vari\'et\'e projective int\`egre $X\subset \P^N$  au
point $x$ et  plus g\'en\'eralement, $T^{k}_{x}X\subset
\P^N$ son espace
$k-$tangent   au point
$x$ (\'etant donn\'ee une repr\'esentation param\'etrique locale au point $x$, c'est l'espace engendr\'e par $x$ et les
d\'eriv\'ees partielles d'ordre inf\'erieur ou \'egal  \`a $k$). La vari\'et\'e duale
$X^{\vee}\subset
\P^{N\vee}$ est d\'efinie comme la cl\^oture de Zariski de l'ensemble 
$\{H\mid T_xX\subset H, X\,\,{\rm lisse}\,\,{\rm en}\,\,x\}$. On d\'efinit de la m\^eme mani\`ere les vari\'et\'es 
duales sup\'erieures $$X^{k\vee}=
\overline{\{H\mid T^{k}_xX\subset H, X\,\,{\rm lisse}\,\,{\rm en}\,\,x\}}$$ 
\section{Surface de Togliatti}
\label{presentationtogliatti}
\begin{defi}
La surface rationnelle, image du plan projectif par  l'application d\'ecrite ci-dessous 
$$\P^{2}\longrightarrow \P^{6}, (x_{0},x_{1},x_{2}) \mapsto 
(x_0x_1x_2, x_0^2x_1,x_0^2x_2,x_1^2x_2,x_1x_2^2,x_0x_1^2,x_0x_2^2)$$ 
est appel\'ee Surface de Togliatti.
\end{defi}
Les sept formes cubiques du vecteur image s'annulant simultan\'ement aux points $$e_{1}=(1,0,0), e_{2}=(0,1,0), e_{3}=(0,0,1),$$ cette surface est la surface de Del Pezzo, not\'ee $S_{6}$, obtenue en \'eclatant $\P^{2}$ le long de ces trois points par le syst\`eme lin\'eaire des cubiques. L'hyperplan osculateur en un point $P=(x_0,x_1,x_2)$  s'interpr\`ete alors comme une cubique passant  par les points $e_{1},e_{2},e_{3}$ et triple au point $P$, plus pr\'ecis\'ement comme la r\'eunion  des trois droites joignant $P$ aux trois points base; une  \'equation  \'etant 
$$ (x_{2}X_{1}-x_{1}X_{2})(x_{2}X_{0}-x_{0}X_{2})(x_{1}X_{0}-x_{0}X_{1})=0.$$
Comme l'ont remarqu\'e Lanteri et Mallavibarena (voir \cite{LM}, pages 357-359), cette \'equation d\'evelopp\'ee  ne s'exprime qu'avec les $6$ formes cubiques suivantes
$$X_0^2X_1,X_0^2X_2,X_1^2X_2,X_1X_2^2,X_0X_1^2,X_0X_2^2.$$ Notons $V$ l'espace vectoriel engendr\'e par ces six formes cubiques. Alors,  dans 
l'espace projectif $\P^{6}=\P((X_0X_1X_2 \oplus V)^{\vee})$ des cubiques s'annulant aux points $e_{1},e_{2},e_{3}$, les hyperplans 
osculateurs de $S_{6}$ passent par le point $(1,0,0,0,0,0)$. 

\medskip

Nous venons de rappeler comment est d\'efinie la surface de Togliatti ainsi que la propri\'et\'e d'incidence de ses hyperplans osculateurs. Cette surface peut aussi \^etre d\'ecrite comme une section hyperplane g\'en\'erale de $Seg(1,3)\subset \P^{7}$ qui est l'image par le morphisme de Segre de $\P^{1}\times \P^{1}\times \P^{1}$.
 \begin{pro}
Soit $H$ un hyperplan  g\'en\'eral de $\P^7$. Alors $H\cap Seg(1,3)\simeq S_{6}$. 
\end{pro}
{\bf Preuve.} Sur $\P^1\times \P^1\times \P^1$ l'\'equation de  $H$ est de la forme suivante 
$$\phi((X_{i}Y_{j}Z_{k})_{0\le i,j,k\le 1})
=(\sum a_{i,j}Y_iZ_j)X_0+(\sum b_{i,j}Y_iZ_j)X_1$$
Comme les biformes $\sum a_{i,j}Y_iZ_j$ et $\sum b_{i,j}Y_iZ_j$ s'annulent  simultan\'ement en deux points, la surface $H\cap Seg(1,3)$ est  l'\'eclatement de $\P^1\times \P^1$ le long de deux points, i.e. l'\'eclatement de $\P^2$ en trois points (non
 align\'es).~$\Box$
 
\medskip

\textbf{Remarque :}
Le th\'eor\`eme de bidualit\'e  (voir par exemple \cite{Ha} thm 15.24) ne s'\'etend pas trivialement aux ordres  sup\'erieurs de tangence. 
La surface de Togliatti, par exemple, ne v\'erifie pas un  r\'esultat attendu  de ``biosculation'' car la surface de ses hyperplans osculateurs est d\'eg\'en\'er\'ee. Cette derni\`ere   est, \`a ma connaissance le
seul exemple un peu \'elabor\'e de surface de $\P^{6}$ telle que  $(S^{3\vee})^{3\vee}\neq S$.

\medskip

Le \textbf{th\'eor\`eme de Togliatti} classifie les r\'eseaux de cubiques 
$\bigwedge \subset H^{0}O_{\P^{2}}(3) $ pour lesquels la surface  image de  la surface de Veronese $v_{3}(\P^{2})$ via la projection  $\P^{9}\setminus \bigwedge \rightarrow \P^{6}$ v\'erifie une \'equation de Laplace (i.e. ses hyperplans osculateurs ont un point commun). 
Deux cas se pr\'esentent : soit le r\'eseau poss\`ede une conique fixe soit non. Le premier cas est r\'esolu ais\'ement. Le second n'est pas du tout trivial et  Togliatti montre que 
la projection de  $v_{3}(\P^{2})$ v\'erifie une \'equation de Laplace si et seulement si ce r\'eseau est engendr\'e par les cubes de trois formes lin\'eaires (voir \cite{FI} et \cite{T}).

\smallskip 

Le th\'eor\`eme ci-dessous est une formulation simplifi\'ee du th\'eor\`eme de Togliatti. 
Seul le cas du r\'eseau non tangent y est trait\'e et on suppose m\^eme que le $\P^{3}$ contenant le r\'eseau (qui se projette sur le point commun des hyperplans osculateurs) n'a pas de lieu base (Pour une exploration exhaustive voir le papier 
\cite{FI} de D.Franco et G. Ilardi).  D\'eterminer les projections de $v_{3}(\P^{2})$ qui v\'erifient une \'equation de Laplace revient alors \`a classifier les  $\P^{3}\subset \P(H^{0}O_{\P^{2}}(3))$ qui coupent la vari\'et\'e des cubiques d\'ecompos\'ees le long d'une surface.
\begin{theo}\label{toto}
Soient $f_i\in \C[X_{0},X_{1},X_{2}], i=0,\cdots, 3$ quatre polyn\^omes homog\`enes de degr\'e 3 sans z\'ero commun. Supposons que pour toute droite  d'\'equation $l=0$ il existe $(\alpha_0,\alpha_1,\alpha_2,\alpha_3)\in \C^{4}\setminus \lbrace 0\rbrace$ tel que $l$ divise $\sum \alpha_if_i$. Alors 
il existe des formes lin\'eaires $l_0,l_1,l_2$ telles que 
$$\mathrm{vect}(f_0,f_1,f_2,f_3)=\mathrm{vect}(l_0^{3},l_1^{3},l_2^{3},l_{0}l_{1}l_{2})$$
\end{theo} 
{\bf Preuve.} 
Consid\'erons le morphisme $\P^2 \longrightarrow \P^3$ donn\'ee par les
quatre formes
cubiques. Comme pour toute droite $L\subset  \P^2$ d'\'equation $l=0$  il existe une cubique du
syst\`eme
lin\'eaire d'\'equation $lq=0$  o\`u $q=0$ est l'\'equation d'une conique $C\subset \P^{2}$ on en d\'eduit
que l'image de
$L$, ainsi que l'image de $C$ est une cubique plane donc une cubique singuli\`ere. Par cons\'equent
 sur toute
 droite deux points ont la m\^eme image. Le degr\'e du morphisme est donc
 au moins \'egal
 \`a $2$ (un syst\`eme lin\'eaire ne contractant pas de courbe et le
 cardinal d'une fibre
 est minimal sur un ouvert). Comme ce degr\'e doit diviser
 $9$ on en d\'eduit qu'il est
 \'egal
 \`a $3$ ($9$ est exclu sinon nous avons un r\'eseau de cubiques passant
 par  $9$
 points). 
 
 \smallskip
 
 Montrons que la fibre d'un point g\'en\'eral est form\'ee de trois points non
 align\'es. 
 
 \smallskip
 
S'ils \'etaient align\'es, quitte \`a changer de coordonn\'ees nous aurions 
$$ \mathrm{vect}(f_0,f_1,f_2,f_3)=\mathrm{vect}(x_{0}q,g_{1},g_{2},g_{3})$$ avec $q$ de degr\'e 2, $g_{i}$ de degr\'e 3 et 
 $g_{1}(0,x_{1},x_{2})=g_{2}(0,x_{1},x_{2})=g_{3}(0,x_{1},x_{2})$. Ceci implique 
 $g_{i}(x_{0},x_{1},x_{2})=x_{0}q_{i}(x_{0},x_{1},x_{2})+g(x_{0},x_{1},x_{2})$
 avec $g$ de degr\'e 3, $q_{i}$ de degr\'e 2.
 Mais alors les quatre formes $f_0,f_1,f_2,f_3$ ont un z\'ero commun sur $x_{0}=0$.
 
 \smallskip

Les trois
droites joignant les trois points de la fibre ont alors pour image la cubique
singuli\`ere. On en d\'eduit que $C$ est form\'ee des deux autres droites.
Comme l'application  de $L$ sur la cubique singuli\`ere image $\Gamma$ est
birationnelle
l'application de la conique r\'esiduelle $C$ sur
$\Gamma $ est un
rev\^etement double (i.e. \`a chaque point de $l$ on associe deux points
sur $C$).

\smallskip

Si le
rev\^etement triple $\P^2 \rightarrow S $, o\`u $S$ est la cubique image,
 est galoisien son
 groupe de Galois est
 $\Z/3\Z$, dans une base adapt\'ee, en tant qu'automorphisme de $\P^2$,
 il est diagonalisable et est engendr\'e par
 $$g=\left (
 \begin{array}{ccc}
 1      &  0 & 0\\
 0     & j &0\\
 0     &  0  &j^2
 \end{array} \right )
 $$
 Comme $\P^3={\rm Proj} \C[X_0^3,X_1^3,X_2^3,X_0X_1X_2]$ est le seul $\P^3$
 de cubiques sur
 lequel le groupe ci-dessus agisse trivialement le r\'esultat est prouv\'e
 lorsque
 l'extension est galoisienne.

 \smallskip

 Si l'extension n'est pas galoisienne, on consid\`ere la cl\^oture
 galoisienne $\pi : X
 \rightarrow \P^2$.
 C'est un rev\^etement double de $\P^2$ ramifi\'e et un rev\^etement  de
 degr\'e
 $6$ de la cubique image $S$,
 $$ \begin{CD}
 X @>{\pi}>> \P^2\\
 @VVV @VpVV \\
 \pi^{-1}(S)@>{\pi}>>S
 \end{CD}
 $$
 Comme $X$ est normale l'image inverse
 $\pi^{-1}(L)$ d'une droite g\'en\'erale de $\P^2$ est irr\'eductible
 d'apr\`es le th\'eor\`eme de Bertini (thm 8.18 page 179 \cite{Har}). On
 consid\`ere la  transposition $\sigma$ qui \'echange les deux points de la
 fibre par
 $\pi$. L'image par $\pi$ de  la courbe  
 $\sigma(\pi^{-1}(L))$ est une courbe  irr\'eductible  telle que 
 $\pi(\sigma(\pi^{-1}(L))\cup \pi^{-1}(L) )=p^{-1}(p(L))=L\cup C$, il s'agit donc de $C$.
 En effet au
 dessus de chaque point de $L$ la fibre est constitu\'ee des deux points
 tels que le
 triplet est la fibre par $p$. Mais ceci contredit le fait que
 $C$ est form\'ee de deux droites.~$\Box$
\begin{coro}
Soit $E$ le fibr\'e vectoriel de la suite exacte ci-dessous
$$ \begin{CD}
 0 @>>> E  @>>>
 O_{\P^{2\vee}}^{4}  @>>> O_{\P^{2\vee}}(3)
 @>>>0
 \end{CD}
 $$ 
Soit $l$ une droite g\'en\'erale de $\P^{2}$. Alors 
 $E_{\mid l}={\cal O}_{l}\oplus {\cal O}_{l}(-1)\oplus {\cal O}_{l}(-2)$
si et seulement s'il existe trois formes lin\'eaires $l_{0},l_{1},l_{2}$ telles que 
$$ \begin{CD}
 0 @>>> E  @>>>
 O_{\P^{2\vee}}^{4}  @>(l_0^{3},l_1^{3},l_2^{3},l_{0}l_{1}l_{2})>> O_{\P^{2\vee}}(3)
 @>>>0
 \end{CD}
 $$
\end{coro}
\section{Polarit\'e par rapport \`a une courbe rationnelle normale}
\begin{theo}\label{vero}
Soient $2n+1$ points de $v_{2n+1}(\P^{2})$ en position g\'en\'erale, $\frak{P}=\P^{2n}$ l'espace projectif qu'ils engendrent et  $S_{2n(2n+1)}$ l'image de  $v_{2n+1}(\P^{2})$  par la projection de centre $\frak{P}$. Sous ces hypoth\`eses, 
les hyperplans $2n$-tangents de $S_{2n(2n+1)}$ ont un point commun. 
\end{theo}
\textbf{Preuve.} Soient $l_{0}, \cdots, l_{2n+1}$ des formes lin\'eaires en position g\'en\'erales. Afin de montrer que les hyperplans $2n$-tangents de $S_{2n(2n+1)}$ ont un point commun on montre que  l'espace projectif 
$\P=\P(l_{0}^{2n+1}, \cdots, l_{2n+1}^{2n+1}, \prod_{i}\overline{l_{i}})$ rencontre tous les espaces $2n$-tangents de $v_{2n+1}(\P^{2})$. L'image de $\P$ est le point commun recherch\'e.

\smallskip

Soit $l$ une forme lin\'eaire sur $\P^{2}$ et $L$ la droite d'\'equation $l=0$.
En un point $l^{2n+1}$ de la Veronese l'espace $2n$-tangent est
donn\'e par
l'ensemble des formes de degr\'e
$2n+1$ qui
sont divisibles par
$l$  i.e. par les formes de degr\'e $2n$. Notons $U\simeq \C^{2}$ un espace vectoriel tel que  $L=\P U$ et $C_{2n+1}\subset \P S^{2n+1}U$ l'image de $L$ par le plongement de Veronese. On note $\overline{l}_{i}$ la restriction des formes lin\'eaires modulo $l$. 
D'apr\`es le r\'esultat bien connu de polarit\'e des courbes rationnelles normales
(qui affirme que  le point d'intersection des hyperplans $2n$-tangents de $C_{2n+1}$ en $(2n+1)$ points g\'en\'eraux distincts appartient au $\P^{2n}$ engendr\'e par ces $(2n+1)$ points)
le sous espace projectif de dimension $2n$ de $\P S^{2n+1}U$  engendr\'e par les $(2n+1)$ points 
$(\overline{l}_{0}^{2n+1}, \cdots,\overline{l}_{2n+1}^{2n+1}) $ contient le point
$\prod_{i}\overline{l_{i}}$ d'intersection des $(2n+1)$ hyperplans $2n$-tangents de $C_{2n+1}$. Modulo $l$ les $2n+1$ formes  sont donc li\'ees, ce qui prouve le th\'eor\`eme.~$\Box$
\section{Polarit\'e par rapport \`a un produit de droites projectives}
Soit $U$ un espace vectoriel complexe de dimension $2$. 
Notons 
$Seg(1,n)\subset \P U^{\otimes n}$ l'image de $\P U\times \cdots \times\P U$ par le plongement de Segre.
\begin{lem}
$Seg(1,n)^{(n-1)\vee}\simeq Seg(1,n)$.
\end{lem}
{\bf Preuve.} Pour s'en convaincre donnons-en une preuve d\'etaill\'ee pour $Seg(1,3)$.
Notons 
$$\P^{7}=
\P(\C[X_{0}Y_{0}Z_{0},X_{0}Y_{0}Z_{1},X_{0}Y_{1}Z_{0},X_{0}Y_{1}Z_{1},X_{1}Y_{0}Z_{0},
X_{1}Y_{1}Z_{0},X_{1}Y_{1}Z_{1}])$$
L'hyperplan d'\'equation $\phi(X_{i}Y_{j}Z_{k})=\sum \alpha_{i,j,k}X_{i}Y_{j}Z_{k}=0$ est osculateur au point 
$(x_{i}y_{j}z_{k})$ si et seulement si 
$$ \frac{\partial^{2}\phi}{\partial X_{i}\partial Y_{j}}(x_{i}y_{j}z_{k})=0,
\frac{\partial^{2}\phi}{\partial X_{i}\partial Z_{k}}(x_{i}y_{j}z_{k})=0,
\frac{\partial^{2}\phi}{\partial Y_{j}\partial Z_{k}}(x_{i}y_{j}z_{k})=0
$$
On obtient imm\'ediatement 
$\phi(X_{i}Y_{j}Z_{k})=(x_{1}X_{0}-x_{0}X_{1})
(y_{1}Y_{0}-y_{0}Y_{1})(z_{1}Z_{0}-z_{0}Z_{1})$.

\smallskip

Plus g\'en\'eralement, 
notons $ \underline X_{i}$  (resp. $ \underline x_{i}$) les coordonn\'ees 
$X_{i,0}, X_{i,1}$ (resp. les nombres $x_{i,0},x_{i,1}$ ). 
L'espace $n$-tangent en un point de $Seg(1,n)$ est donn\'e par les d\'eriv\'ees $(n-1)-i\grave{e}mes$ d'un syst\`eme de coordonn\'ees locales qui est, essentiellement,  le produit des variables. Ainsi
la forme multilin\'eaire 
  $\phi(\underline X_{1}, \underline X_{2},\cdots, \underline X_{n})$
est la trace sur  $Seg(1,n)$ d'un hyperplan $n$-tangent  
$ H\subset \P^{2^{n}-1}$ au point $(\underline x_{1}, \underline x_{2},\cdots , 
\underline x_{n})$ si et seulement si 
$$ \phi(\underline x_{1}, \underline X_{2},\cdots, \underline X_{n})=\phi(\underline X_{1}, \underline x_{2},\cdots , \underline X_{n})=\cdots =
 \phi(\underline X_{1}, \underline X_{2},\cdots,\underline X_{n-1},\underline x_{n})=0$$
Apr\`es un calcul \'el\'ementaire on en d\'eduit
$$\phi(\underline X_{1}, \underline X_{2},\cdots, \underline X_{n})=\prod_{i=1,\cdots ,n} (x_{i,1}X_{i,0}-x_{i,0}X_{i,1})$$
Ce qui prouve le lemme.~$\Box$

\medskip

La polarit\'e point-hyperplan par rapport \`a $C_{n}$ permet d'\'etendre 
l'isomorphisme $C_{n}^{n\vee}\simeq C_n$ au $\P^{n}$ ambiant. D'une fa\c{c}on similaire (une sorte de 'polarit\'e' par rapport au Segre) l'isomorphisme  $Seg(1,n)^{(n-1)\vee}\simeq Seg(1,n)$ s'\'etend aux espaces projectifs ambiants. Ici aussi il faudra distinguer le cas pair du cas impair. 
\begin{pro} 
\label{prod} L'isomorphisme  $Seg(1,n)\simeq Seg(1,n)^{(n-1)\vee}$ qui \`a un point du Segre associe l'hyperplan $n$-tangent au Segre en ce point s\'etend aux espace projectifs 
$\P(U^{\otimes n})$ et $\P(U^{*\otimes n})$ de la mani\`ere suivante : Au point 
$x\in \P(U^{\otimes n})$ on associe le point d'intersection des hyperplans  
$n$-tangents aux points de $x^{\vee}\cap Seg(1,n)^{(n-1)\vee}$.
Lorsque $n$ est impair, ce point image appartient \`a $x^{\vee}$. Lorsque $n$ est pair, les points $x\in \P(U^{\otimes n})$ pour lesquels le point image appartient \`a $x^{\vee}$ forment une hyperquadrique.   
\end{pro}
{\bf Preuve.} 
L'isomorphisme  $m^{\otimes n} : U^{\otimes n}\rightarrow U^{*\otimes n}$
donn\'e par la matrice $m^{\otimes n}=\left (
\begin{array}{cc}
 0& -1\\
 1& 0
\end{array}
\right)^{\otimes n} $ est sym\'etrique pour $n$ pair et antisym\'etrique pour $n$ impair.  Pour $x\in \P(U^{\otimes n})$ on note $x^{\vee}$ l'\'equation de l'hyperplan correpondant dans $ \P(U^{*\otimes n})$.  
Notons $\xi \in x^{\vee}\cap Seg(1,n)^{(n-1)\vee}$ le point g\'en\'erique de l'intersection. On a par d\'efinition 
$<x,\xi>=0$, o\`u $<.,.>$ est le crochet de dualit\'e. Comme $m^{\otimes n}$ est inversible il existe un unique 
$x_{\xi}\in Seg(1,n)$ tel que $\xi=m^{\otimes n}(x_{\xi})$.  On en d\'eduit,  
en utilisant les propri\'et\'es de sym\'etrie de la matrice $m^{\otimes n}$, que 
$<x,m^{\otimes n}(x_{\xi})>=<x_{\xi},m^{\otimes n}(x)>=0$,  c'est \`a dire  que le point $m^{\otimes n}(x)$ appartient \`a l'hyperplan $n$-tangent g\'en\'erique de la vari\'et\'e 
$x^{\vee}\cap Seg(1,n)^{(n-1)\vee}$. 

\smallskip
 
Enfin lorsque $n$ est impair, la matrice $m^{\otimes n}$ \'etant antisym\'etrique, on a 
$<x,m^{\otimes n}(x)>=0$ pour tout $x\in \P(U^{\otimes n})$. Lorsque $n$ est pair la matrice est sym\'etrique et les point $x\in \P(U^{\otimes n})$ tels que $<x,m^{\otimes n}(x)>=0$ forment une quadrique de $\P(U^{\otimes n})$.~$\Box$
\begin{theo}\label{pun}
Soit $X_{2n}$ une section hyperplane g\'en\'erale de $Seg(1,2n+1)$. Les hyperplans $2n$-tangents de  $X_{2n}$  ont un point commun.
\end{theo}
\textbf{Preuve.} Les hyperplans $2n$-tangents de $X_{2n}$ sont les hyperplans tangents de $Seg(1,2n+1)$ coup\'es par l'hyperplan consid\'er\'e. Le th\'eor\`eme est alors  une cons\'equence imm\'ediate de la proposition \ref{prod}.~$\Box$


\medskip

Jean Vall\`es\\
Laboratoire de Math\'ematiques appliqu\'ees\\ 
de Pau et des Pays de l'Adour\\
Avenue de l'universit\'e\\
64000 Pau, France\\
jean.valles@univ-pau.fr
compil\'e le 3 avril  2006

\end{document}